\documentclass{amsart}
\pdfoutput=1
\usepackage[a4paper, margin=2cm]{geometry}
\usepackage[utf8]{inputenc}
\usepackage{amsthm,thmtools,amssymb}
\usepackage{bussproofs}
\usepackage{enumitem}
\usepackage[pagebackref]{hyperref}
\usepackage{comment}
\hypersetup{
    colorlinks=true,
    allcolors={purple}
}
\usepackage{tikz}
\usetikzlibrary{arrows.meta,decorations.pathreplacing,decorations.markings,shapes,calc}
\usetikzlibrary{matrix,arrows,decorations.pathmorphing,backgrounds,decorations.markings,positioning}
\tikzset{2cell/.style={-implies,double,double equal sign distance,shorten 
>=2pt, shorten <=3pt}}
\tikzset{2cellshort/.style={-implies,double,double equal sign distance,shorten 
>=4pt, shorten <=5pt}}
\tikzset{2cellr/.style={implies-,double,double equal sign distance,shorten 
>=3pt, shorten <=2pt}}
\tikzset{3cell/.style={-implies,double,double distance=2.5pt,shorten >=2pt, 
shorten <=3pt}}
\tikzset{labelsize/.style={font=\scriptsize}}
\tikzset{string/.style={very thick}}
\tikzset{
  pto/.style={->,postaction={decorate},
    decoration={
        markings,
        mark=at position 0.5 with {\arrow{|}}}
  },
}

\mathchardef\mhyphen="2D

\setlist[description]{font=\normalfont}

\declaretheorem[style=plain,numberwithin=section,name=Theorem]{theorem}
\declaretheorem[style=plain,sibling=theorem,name=Lemma]{lemma}
\declaretheorem[style=plain,sibling=theorem,name=Proposition]{proposition}
\declaretheorem[style=plain,sibling=theorem,name=Corollary]{corollary}

\declaretheorem[style=definition,sibling=theorem,name=Example]{example}

\declaretheorem[style=definition,sibling=theorem,name=Remark]{remark}

\makeatletter
\newcommand{\pto}{}
\newcommand{\pgets}{}

\DeclareRobustCommand{\pto}{\mathrel{\mathpalette\p@to@gets\to}}
\DeclareRobustCommand{\pgets}{\mathrel{\mathpalette\p@to@gets\gets}}

\newcommand{\p@to@gets}[2]{%
  \ooalign{\hidewidth$\m@th#1\mapstochar\mkern5mu$\hidewidth\cr$\m@th#1\to$\cr}%
}
\makeatother

\newcommand{\cat}[1]{\mathcal{#1}}

\newcommand{\defemph}[1]{\textbf{#1}}

\newcommand{\op}{\mathrm{op}}

\newcommand{\OSRsa}{\mathbf{OSR_{\mathrm{sa}}}}

\newcommand{\DLat}{\mathbf{DLat}}

\newcommand{\Frm}{\mathbf{Frm}}

\newcommand{\IQuant}{\mathbf{IQuant}}

\newcommand{\Idl}{\mathrm{Idl}}
\newcommand{\Rad}{\mathrm{Rad}}

\newcommand{\Rnonneg}{{\mathbb{R}}_{\geq 0}}
\newcommand{\Rbar}{\overline{\mathbb{R}}}
\newcommand{\Spec}{\mathrm{Spec}}
\newcommand{\pt}{\mathrm{pt}}

\title{Ordered semirings and subadditive morphisms}
\author{Soichiro Fujii}
\thanks{The author gratefully acknowledges the support of JSPS Overseas Research Fellowships}
\address{School of Mathematical and Physical Sciences, Macquarie University, NSW 2109, Australia}
\email{s.fujii.math@gmail.com}
\date{November 7, 2023}

\begin{document}

\begin{abstract}
An ordered semiring is a commutative semiring equipped with a compatible preorder. Ordered semirings generalise both distributive lattices and commutative rings, and provide a convenient framework to unify certain aspects of lattice theory and ring theory. The ideals of an ordered semiring $A$ form a commutative integral quantale $\mathrm{Idl}(A)$, and similarly, the radical ideals of $A$ form a (spatial) frame $\mathrm{Rad}(A)$. We characterise $\mathrm{Idl}$ and $\mathrm{Rad}$ as the left adjoints of the (non-full) inclusion functors from the categories of commutative integral quantales and of frames, respectively, to that of ordered semirings and subadditive morphisms between them. The (sober) topological space $\mathrm{pt}(\mathrm{Rad}(A))$ corresponding to $\mathrm{Rad}(A)$ is homeomorphic to the space $\mathrm{Spec}(A)$ of prime ideals of $A$.
\end{abstract}

\maketitle

\section{Introduction}
The notion of \emph{ideal}, as well as its prime and maximal variants, appears both in lattice theory and in ring theory.
What is the relationship between lattice-theoretic ideals and ring-theoretic ideals?
A common answer to this question is that they coincide in the special cases of \emph{Boolean algebras} and \emph{Boolean rings}---it is well-known that Boolean algebras and Boolean rings are equivalent objects (see e.g.~\cite[Exercise 1.24]{Atiyah-Macdonald} or \cite[Theorem~I.1.9]{Johnstone_Stone_spaces}). 
However, this does not fully explain why these notions continue to behave so similarly beyond these classes.
The purpose of this short note is to give a conceptual explanation of this, by presenting a unifying framework.

We point out that both distributive lattices and rings are special cases of \emph{ordered semirings}, (commutative) semirings equipped with a compatible preorder.
Ordered semirings, when considered with \emph{subadditive morphisms} rather than the usual (monotone) homomorphisms as morphisms between them, provide a convenient framework to develop a general theory of ideals.

We define ideals of an ordered semiring, generalising ideals of a distributive lattice and of a ring.
For any ordered semiring $A$, the poset $\Idl(A)$ of its ideals (ordered by inclusion) is a complete lattice.
One can moreover take the \emph{product} $I\cdot J$ of two ideals $I,J\in\Idl(A)$, as in the case of ideals of a ring.
With these operations, $\Idl(A)$ forms a structure known as a (commutative) \emph{integral quantale}.
Indeed, we obtain a functor $\Idl\colon \OSRsa\to \IQuant$, from the category $\OSRsa$ of ordered semirings and subadditive morphisms between them, to the category $\IQuant$ of integral quantales and quantale homomorphisms between them.
On the other hand, integral quantales are a special case of ordered semirings, and hence we also have the (non-full) inclusion functor $\IQuant\to \OSRsa$.
It turns out that $\Idl$ is the left adjoint of this inclusion functor.
Notice that since rings do not contain integral quantales as a special case, such a simple characterisation of $\Idl$ is not available in the classical ring-theoretic framework.

Integral quantales in turn include \emph{frames}, which are certain structures used in a point-free approach to topology (see e.g.~\cite[Chapter~II]{Johnstone_Stone_spaces}).
The inclusion functor $\Frm\to \OSRsa$ from the category $\Frm$ of frames and frame homomorphisms also has a left adjoint; this is a formal consequence of the above observation and the known result that the inclusion $\Frm \to \IQuant$ has a left adjoint \cite{Niefield_Rosenthal}.
Explicitly, the left adjoint of the inclusion functor $\Frm\to\OSRsa$ is the functor $\Rad\colon \OSRsa\to\Frm$ which maps each ordered semiring $A$ to the frame $\Rad(A)$ of \emph{radical ideals} of $A$.
Frames of the form $\Rad(A)$ are \emph{spatial}, meaning that they correspond to certain (sober) topological spaces. 
In fact, this generalises known approaches to defining spectra of distributive lattices (for which $\Rad(A)=\Idl(A)$ holds) or of rings (see e.g.~\cite[Section~II.3]{Johnstone_Stone_spaces} or \cite[Section~V.3]{Johnstone_Stone_spaces}, respectively).

We finally show that the topological space $\pt(\Rad(A))$ corresponding to the frame $\Rad(A)$ of radical ideals of an ordered semiring $A$ has another description in terms of prime ideals. 
As in the case of distributive lattices, we can identify the prime ideals of an ordered semiring $A$ with the subadditive morphisms from $A$ to a certain two-element ordered semiring $\mathbf{2}$.
We can thus regard the set of all subadditive morphisms of type $A\to \mathbf{2}$ as the subspace of the product space $\mathbf{2}^A$, where $\mathbf{2}$ is now being regarded as the Sierpi\'nski space. This defines the topological space $\Spec(A)$ of prime ideals of $A$, which is canonically homeomorphic to $\pt(\Rad(A))$.

\subsection*{Related work}
Modulo minor differences, ordered semirings are called \emph{halos} in \cite{Paugam}, where subadditive morphisms and more general morphisms are considered. \cite{Paugam} also discusses localisations of ordered semirings and many examples.

\subsection*{Conventions}
All lattices are assumed to be bounded, and all lattice homomorphisms preserve the bottom and top elements.
All rings and semirings are assumed to be unital and commutative, and all ring and semiring homomorphisms preserve the multiplicative unit. 

\section{Ordered semirings}
Recall that a (commutative) \defemph{semiring} is a tuple $(A,0,1,+,\cdot)$ where $(A,0,+)$ and $(A,1,\cdot)$ are commutative monoids, such that 
\[
x\cdot 0=0
\qquad \text{and} \qquad 
x\cdot (y+z)=x\cdot y+x\cdot z
\]
hold for all $x,y,z\in A$.
An \defemph{ordered semiring} is a tuple $(A,\leq,0,1,+,\cdot)$ where $(A,\leq)$ is a preordered set\footnote{All results in this paper continue to hold even if we required $(A,\leq)$ to be a poset.} and $(A,0,1,+,\cdot)$ is a semiring, such that $+$ and $\cdot$ are monotone with respect to $\leq$. 

\begin{example}
\label{ex:semiring}
    Any semiring can be regarded as a discretely ordered (i.e.~$x\leq y$ only if $x=y$) semiring. In particular, all rings and fields induce discretely ordered semirings.
\end{example}
\begin{example}
\label{ex:dist_lat}
    Any distributive lattice $(D,\bot,\top,\vee,\wedge)$ can be regarded as an ordered semiring $(D,{\leq},\bot,\top,\vee,\wedge)$, where $\leq$ is the partial order on $D$ with respect to which $\vee$ is the join operation (and $\wedge$ the meet).
    In particular, the two-element chain $\mathbf{2}=(\{0,1\}, \leq, 0,1, \vee,\wedge)$ is an ordered semiring.

    Recall that a \defemph{frame} is a complete lattice $(F,\leq)$ such that 
    \[
x\wedge\Big(\bigvee_{\lambda\in\Lambda}y_\lambda\Big)=\bigvee_{\lambda\in\Lambda}x\wedge y_\lambda
    \] 
    holds for each set $\Lambda$, $x\in F$, and $\Lambda$-indexed family $(y_\lambda)_{\lambda\in\Lambda}$ of elements of $F$. Since every frame $(F,\leq)$ is a distributive lattice, it induces an ordered semiring $(F,{\leq},\bot,\top,\vee,\wedge)$.
\end{example}
\begin{example}
\label{ex:idempotent-semiring}
    Any semiring $(A,0,1,+,\cdot)$ in which $x+x=x$ holds for all $x\in A$ induces an ordered semiring $(A,\leq,0,1,+,\cdot)$, where we define $x\leq y\iff x+y=y$. Notice that $(A,0,+)$ is an idempotent commutative monoid, i.e., a semilattice. The addition is the join operation in the poset $(A,\leq)$.

    Recall that a (unital and commutative) \defemph{quantale} is a tuple $(Q,\leq,1,\cdot)$ where $(Q,\leq)$ is a complete lattice and $(Q,1,\cdot)$ is a commutative monoid, such that 
    \[
    x\cdot\Big(\bigvee_{\lambda\in\Lambda}y_\lambda\Big)=\bigvee_{\lambda\in\Lambda}x\cdot y_\lambda
    \]
    holds for each set $\Lambda$, $x\in Q$ and $\Lambda$-indexed family $(y_\lambda)_{\lambda\in\Lambda}$ of elements of $Q$. 
    A quantale $(Q,\leq,1,\cdot)$ is called \defemph{integral} (or affine) if $1$ is the largest element of $(Q,\leq)$. Since any quantale $(Q,\leq,1,\cdot)$ induces a semiring $(Q,\bot,1,\vee,\cdot)$ in which addition is idempotent, it also induces an ordered semiring $(Q,\leq,\bot,1,\vee,\cdot)$.
\end{example}
\begin{example}
    Nonnegative real numbers $\Rnonneg=([0,\infty),\leq,0,1,+,\cdot)$ and natural numbers $\mathbb{N}=(\mathbb{N},\leq,0,1,+,\cdot)$ with the usual order, addition and multiplication, form ordered semirings. 
    Notice that $(\mathbb{R},\leq,0,1,+,\cdot)$ is not an ordered semiring because $\cdot$ is not monotone for negative real numbers.

    More generally, consider the tuple $(R,{\leq}, 0,1,-,+,\cdot)$ where $(R,\leq)$ is a preordered set and $(R,0,1,-,+,\cdot)$ is a ring such that $+$ is monotone, $0\leq 1$, and $0\leq x$ and $0\leq y$ imply $0\leq x y$. Then the nonnegative cone $R_{\geq 0}=\{\,x\in R\mid 0\leq x\,\}$ of $R$ forms
    an ordered semiring. 
\end{example}
\begin{example}
\label{ex:Gamma}
    Real numbers augmented with the least element $-\infty$ form the \defemph{max-plus semiring} $\Rbar=(\mathbb{R}\cup\{-\infty\},\leq,-\infty,0,\max,+)$, which is isomorphic to the ordered semiring $([0,\infty), \leq,0,1,\max,\cdot)$.
    
    More generally, let $\Gamma=(\Gamma,\leq,1,\cdot)$ be a \defemph{totally ordered commutative monoid}, i.e., $(\Gamma,\leq)$ is a totally ordered set and $(\Gamma, 1,\cdot)$ is a commutative monoid such that $\cdot$ is monotone. Then the totally ordered set $(\Gamma\cup\{-\infty\},\leq)$ obtained from $(\Gamma,\leq)$ by adjoining the least element $-\infty$, together with the obvious extension of $\cdot$ (i.e., $x\cdot (-\infty)=-\infty$ for all $x\in\Gamma\cup\{-\infty\}$), give rise to an ordered semiring $\overline{\Gamma}=(\Gamma\cup\{-\infty\},{\leq},-\infty,1,\max,\cdot)$. 
    (In fact, this is a special case of Example~\ref{ex:idempotent-semiring}.)
\end{example}

\section{Subadditive morphisms}
Given ordered semirings $A=(A,\leq_A,0_A,1_A,+_A,\cdot_A)$ and $B=(B,\leq_B,0_B,1_B,+_B,\cdot_B)$, a \defemph{subadditive morphism} $f\colon A\to B$ is a monotone function satisfying 
\begin{itemize}
    \item $f(0_A)\leq_B 0_B$;\footnote{All results in this paper continue to hold even if we replaced this axiom by $f(0_A)=0_B$.}
    \item $f(1_A)=1_B$;
    \item for all $x,y\in A$, $f(x+_A y)\leq_B f(x) +_B f(y)$; and
    \item for all $x, y\in A$, $f(x\cdot_A y)=f(x)\cdot_B f(y)$.
\end{itemize}
A \defemph{homomorphism} $f\colon A\to B$ is a subadditive morphism satisfying $f(0_A) = 0_B$ and $f(x+_Ay)=f(x)+_Bf(y)$ for all $x,y\in A$. 
We denote the category of ordered semirings and subadditive morphisms between them by $\OSRsa$.

\begin{remark}
\label{rmk:subadditive-hom-is-hom-when}
A subadditive morphism $f\colon A\to B$ between ordered semirings $A$ and $B$ is necessarily a homomorphism whenever either of the following conditions is satisfied.
\begin{enumerate}
    \item $B$ is discretely ordered.
    \item For each $x\in A$ we have $0_A\leq_A x$, and $(B,0_B,+_B)$ is the join-semilattice for the order $(B,\leq_B)$ (i.e., $B$ is induced from the construction of Example~\ref{ex:idempotent-semiring}). To see this, observe that for all $x,y\in A$ we have $x=x+_A 0_A\leq_A x+_A y$, and hence $f(x)\leq_B f(x+_A y)$; similarly, we have $f(y)\leq_B f(x+_A y)$.
    Thus $f(x)+_Bf(y)=f(x)\vee f(y)\leq_B f(x+_A y)$ holds.
\end{enumerate}
\end{remark}

\begin{example}
By Remark~\ref{rmk:subadditive-hom-is-hom-when} (1), all subadditive morphisms between discretely ordered semirings are homomorphisms, which coincide with semiring homomorphisms. Thus the category of semirings and semiring homomorphisms, as well as the category of rings and ring homomorphisms, are full subcategories of $\OSRsa$ via Example~\ref{ex:semiring}.
\end{example}
\begin{example}
\label{ex:morphism-dist-lat-2}
By Remark~\ref{rmk:subadditive-hom-is-hom-when} (2), all subadditive morphisms between distributive lattices are homomorphisms, which coincide with lattice homomorphisms. Thus the category $\DLat$ of distributive lattices and lattice homomorphisms is a full subcategory of $\OSRsa$ via Example~\ref{ex:dist_lat}. In particular, for a distributive lattice $D$, the subadditive morphisms $D\to \mathbf{2}$ correspond to the (lattice-theoretic) prime ideals of $D$.
\end{example}
\begin{example}
\label{ex:morphism-ring-2}
Let $R$ be a ring, regarded as a discretely ordered semiring. Then the subadditive morphisms $R\to \mathbf{2}$ correspond to the (ring-theoretic) prime ideals of $R$.
\end{example}
\begin{example}
    Let $R$ be an integral domain, and $R[x]$ the polynomial ring over $R$, regarded as a discretely ordered semiring.
    Then the function $\mathrm{deg}\colon R[x]\to \overline{\mathbb{R}}$ mapping each polynomial $p(x)\in R[x]$ to its degree, is a subadditive morphism.
\end{example}
\begin{example}
\label{ex:morphism-Archimedean}
For a field $k$ (regarded as a discretely ordered semiring), the subadditive morphisms $k\to\Rnonneg$ correspond to the absolute values of $k$ \cite[Chapter~XII]{Lang_algebra}.
Similarly, the subadditive morphisms $k\to \Rbar$ correspond to the non-archimedean absolute values of $k$; if we use $\overline{\Gamma}$ for an arbitrary totally ordered abelian group $\Gamma$ instead of $\overline{\mathbb{R}}$, we obtain the valuations of $k$.
See \cite{Paugam} for a further discussion related to this example.
\end{example}

\section{Ideals}
Let $A=(A,\leq,0,1,+,\cdot)$ be an ordered semiring. An \defemph{ideal} of $A$ is a subset $I\subseteq A$ satisfying
\begin{itemize}
    \item for all $x,y\in A$, if $x\leq y$ and $y\in I$, then $x\in I$;
    \item $0\in I$;
    \item for all $x,y\in A$, if $x\in I$ and $y\in I$, then $x+y\in I$; and
    \item for all $x,y\in A$, if $x\in I$ or $y\in I$, then $xy\in I$.  
\end{itemize}
Notice that the ideals of $A$ correspond to the subadditive and submultiplicative morphisms\footnote{For ordered semirings $A$ and $B$, a \defemph{subadditive and submultiplicative morphism} $f\colon A\to B$ is a monotone function satisfying  $f(0_A)\leq_B 0_B$, $f(1_A)\leq_B1_B$, $f(x+_Ay)\leq_B f(x) +_B f(y)$, and $f(x\cdot_A y)\leq_B f(x)\cdot_B f(y)$ for all $x, y\in A$; cf.~\cite{Paugam}.} $A\to\mathbf{2}$.
The set of all ideals of $A$ is denoted by $\Idl(A)$.
The ideals of $A$ are closed under arbitrary intersections and directed unions. For any subset $S\subseteq A$, the smallest ideal $\langle S\rangle$ of $A$ containing $S$ is given by
\[
\langle S\rangle =\{\,z\in A\mid \exists m\in\mathbb{N}.\,\exists s_1,\dots,s_m\in S.\,\exists y_1,\dots,y_m\in A.\,z\leq s_1y_1+\dots+s_my_m\,\}.
\]
For each $x\in A$, the ideal of $A$ generated by $\{x\}$ is denoted by $\langle x\rangle$. Explicitly, we have 
\[
\langle x\rangle =\{\,z\in A\mid \exists y\in A.\,z\leq xy\,\}.
\]
The join operation $\bigvee$ in the complete lattice $(\Idl(A),\subseteq)$ satisfies the following.

\begin{lemma}
    Let $A$ be an ordered semiring, $\Lambda$ a set, and $S_\lambda\subseteq A$ for each $\lambda\in\Lambda$. Then we have 
    \[
        \bigvee_{\lambda\in\Lambda}\langle S_\lambda\rangle = \Big\langle \bigcup_{\lambda\in\Lambda} S_\lambda\Big\rangle.
    \]
\end{lemma}

For $I,J\in\Idl(A)$, we define their \defemph{product} $I\cdot J$ as 
\begin{align*}
I\cdot J&=\bigl\langle\{\,xy\mid x\in I, y\in J\,\}\bigr\rangle \\
&=\{\,z\in A\mid \exists m\in\mathbb{N}.\,\exists x_1,\dots,x_m\in I.\,\exists y_1,\dots,y_m\in J.\, z\leq x_1y_1+\dots+x_my_m\,\}.
\end{align*}

\begin{lemma}
    Let $A$ be an ordered semiring and $S,T\subseteq A$. Then we have $\langle S \rangle \cdot \langle T \rangle =\langle ST\rangle$, where $ST=\{\,st\mid s\in S,t\in T\,\}$.
\end{lemma}
\begin{proof}
    Clearly we have $\langle ST \rangle\subseteq \langle S\rangle \cdot \langle T\rangle$. To show the reverse inclusion, let $x\in \langle S\rangle$ and $y\in\langle T\rangle$. Thus we have 
    \[
        x\leq s_1u_1+\dots + s_mu_m
        \quad\text{and}\quad
        y\leq t_1v_1+\dots + t_nv_n
    \]
    for some $m,n\in\mathbb{N}$, $s_i\in S$, $t_j\in T$ and $u_i,v_j\in A$.
    It follows that 
    \[
        xy\leq \sum_{\substack{1\leq i\leq m \\ 1\leq j\leq n}}s_it_ju_iv_j
    \]
    with $s_it_j\in ST$. Thus we have $xy\in \langle ST\rangle$.
\end{proof} 

It follows that the product of ideals is associative and commutative, and has the largest ideal $A$ as the unit. Moreover, it satisfies the infinitary distributive law
\[
I\cdot \Big(\bigvee_{\lambda\in\Lambda}J_\lambda\Big)=\bigvee_{\lambda\in\Lambda}I\cdot J_\lambda.
\]
Therefore $(\Idl(A),\subseteq,A,\cdot)$ is an integral quantale.
We remark that the idea of using integral quantales as an abstraction of the structure formed by ideals (of a ring) goes back to Krull \cite{Krull-1924}.

It is easy to check that $\langle -\rangle \colon A\to\Idl(A)$ is a subadditive morphism:
\begin{proposition}
    Let $A$ be an ordered semiring and $x,y\in A$. 
    \begin{itemize}
        \item If $x\leq y$, then $\langle x\rangle \subseteq \langle y\rangle$.
        \item $\langle 0\rangle$ is the smallest element of $\Idl(A)$.
        \item $\langle 1\rangle =A$.
        \item $\langle x+y\rangle\subseteq\langle x\rangle \vee\langle y\rangle $.
        \item $\langle xy\rangle =\langle x\rangle \cdot\langle y\rangle$.
    \end{itemize}
\end{proposition}

In fact, it is the universal subadditive morphism from $A$ to an integral quantale:
\begin{theorem}
\label{thm:Idl}
    Let $A$ be an ordered semiring, $Q$ an integral quantale, and $f\colon A\to Q$ a subadditive morphism. Then there exists a unique quantale homomorphism\footnote{A quantale homomorphism is a function between quantales preserving arbitrary joins and the monoid structure.} $g\colon \Idl(A)\to Q$ such that the diagram 
    \begin{equation*}
    \begin{tikzpicture}[baseline=-\the\dimexpr\fontdimen22\textfont2\relax ]
      \node(0) at (0,0) {$A$};
      \node(1) at (3,0) {$\Idl(A)$};
      \node(2) at (3,-1.5) {$Q$};
      \draw [->] (0) to node[auto, labelsize] {$\langle-\rangle$} (1); 
      \draw [->] (0) to node[auto, swap,labelsize] {$f$} (2); 
      \draw [->, dashed] (1) to node[auto,labelsize] {$g$} (2);
    \end{tikzpicture}
    \end{equation*}
    commutes.
\end{theorem}
\begin{proof}
    For each $I\in\Idl(A)$, define 
    \[
    g(I)=\bigvee_{x\in I}f(x).
    \]
    Note that this is the only possible definition to satisfy the stated properties because we have $I=\bigvee_{x\in I}\langle x\rangle$ in $\Idl(A)$.
    
    Observe that for any subset $S\subseteq A$, we have $g\bigl(\langle S\rangle\bigr)=\bigvee_{s\in S}f(s)$. Indeed, for any $z\in \langle S\rangle$, we have $z\leq s_1y_1+\dots+s_my_m$ for some $s_i\in S$ and $y_i\in A$. Thus we have 
    \begin{align*}
    f(z)&\leq f(s_1)f(y_1)\vee\dots \vee f(s_m)f(y_m)\\
    &\leq f(s_1)\vee\dots \vee f(s_m)\\
    &\leq \bigvee_{s\in S}f(s),
    \end{align*}
    where we used the assumption that $Q$ is an integral quantale (which implies $pq\leq p$ for all $p,q\in Q$). 
    
    In particular, for each $x\in A$ we have $f(x)=g\bigl(\langle x \rangle\bigr)$. 
    It is straightforward to check that $g$ is a quantale homomorphism. 
\end{proof} 

Let $\IQuant$ be the category of integral quantales and quantale homomorphisms. Via the construction of Example~\ref{ex:idempotent-semiring}, we have a (non-full) inclusion functor $\IQuant\to \OSRsa$.

\begin{corollary}
    The inclusion functor $\IQuant\to \OSRsa$ has a left adjoint $\Idl\colon\OSRsa\to \IQuant$.
\end{corollary}

\section{Radical ideals}
Let $A$ be an ordered semiring.
A \defemph{radical ideal} of $A$ is an ideal $I\subseteq A$ satisfying
\begin{itemize}
    \item for all $x\in A$ and  $n\in\mathbb{N}$, $x^n\in I$ implies $x\in I$.
\end{itemize}
The set of all radical ideals of $A$ is denoted by $\Rad(A)$.
We shall show that $\Rad(A)$ is a frame, and that moreover $\Rad\colon \OSRsa\to\Frm$ is the left adjoint of the (non-full) inclusion functor $\Frm\to \OSRsa$; here, $\Frm$ is the category of frames and frame homomorphisms\footnote{A frame homomorphism is a function between frames preserving arbitrary joins and finite meets. Note that $\Frm$ is a full subcategory of $\IQuant$.} and the inclusion functor is obtained via the construction of Example~\ref{ex:dist_lat}.

Of course, this can be checked directly as in the previous section. Instead, here we shall recall a result of \cite[Section~3]{Niefield_Rosenthal} constructing the left adjoint of the inclusion functor $\Frm\to \IQuant$ explicitly.
Specifically, if $Q$ is an integral quantale, an element $p\in Q$ is called \defemph{semiprime} if for all $q\in Q$ and $n\in\mathbb{N}$, $q^n\leq p$ implies $q\leq p$. Then the set $S(Q)$ of all semiprime elements forms a frame, and $S\colon \IQuant\to \Frm$ is the left adjoint of the inclusion. (Notice that $S(Q)$ is a subset of $Q$ closed under arbitrary meets. Hence the inclusion $S(Q)\to Q$ has a left adjoint $\sqrt{-}\colon Q\to S(Q)$, which is the component at $Q$ of the unit of the adjunction.)

Now observe that an ideal $I$ of an ordered semiring $A$ is semiprime in $\Idl(A)$ if and only if it is radical. Indeed, every semiprime ideal is radical since $x\in I$ is equivalent to $\langle x\rangle \subseteq I$, and since $\langle x\rangle^n=\langle x^n\rangle$. Conversely, if $I$ is radical and $J$ is an ideal such that $J^n\subseteq I$ for some $n\in\mathbb{N}$, then for each $x\in J$ we have $x^n\in J^n\subseteq I$, and hence $x\in I$. Thus we have $J\subseteq I$.

\begin{theorem}
    Let $A$ be an ordered semiring, $F$ a frame and $f\colon A\to F$ a subadditive morphism. Then there exists a unique frame homomorphism $g\colon \Rad(A)\to F$ such that the diagram 
    \begin{equation*}
    \begin{tikzpicture}[baseline=-\the\dimexpr\fontdimen22\textfont2\relax ]
      \node(0) at (0,0) {$A$};
      \node(1) at (3,0) {$\Rad(A)$};
      \node(2) at (3,-1.5) {$F$};
      \draw [->] (0) to node[auto, labelsize] {$\sqrt{\langle-\rangle}$} (1); 
      \draw [->] (0) to node[auto, swap,labelsize] {$f$} (2); 
      \draw [->, dashed] (1) to node[auto,labelsize] {$g$} (2);
    \end{tikzpicture}
    \end{equation*}
    commutes.
\end{theorem}
\begin{corollary}
    The inclusion functor $\Frm\to \OSRsa$ has a left adjoint $\Rad\colon\OSRsa\to \Frm$.
\end{corollary}

Next we shall determine the class of frames which arise as $\Rad(A)$ for some ordered semiring $A$.
Observe that the inclusion functor $\Frm\to \OSRsa$ factors as 
\[
    \begin{tikzpicture}[baseline=-\the\dimexpr\fontdimen22\textfont2\relax ]
      \node(0) at (0,0) {$\Frm$};
      \node(1) at (2,-1) {$\DLat$};
      \node(2) at (4,0) {$\OSRsa$.};
      \draw [->] (0) to node[auto, labelsize] {} (1); 
      \draw [->] (0) to node[auto, swap,labelsize] {} (2); 
      \draw [->] (1) to node[auto,labelsize] {} (2);
    \end{tikzpicture}
\]
The left adjoint of the inclusion functor $\Frm\to \DLat$ is given by $\Idl$, since every ideal in a distributive lattice is radical. 
Moreover, the inclusion functor $\DLat\to \OSRsa$ also has a left adjoint $L\colon \OSRsa\to \DLat$. Indeed, given any $A\in \OSRsa$, the distributive lattice $LA$ is presented by the generators $\overline{x}$ for each $x\in A$, subject to the relations 
\begin{itemize}
    \item for each $x,y\in A$ with $x\leq y$, $\overline{x}\leq \overline{y}$ (which can be expressed equationally as, for example, $\overline{x} \vee \overline{y}=\overline{y}$);
    \item $\overline{0}=\bot$;
    \item $\overline{1}=\top$;
    \item for each $x,y\in A$, $\overline{x+y}\leq \overline{x}\vee\overline{y}$; and 
    \item for each $x,y\in A$, $\overline{xy} = \overline{x}\wedge \overline{y}$.
\end{itemize}
(See \cite[V.3.1]{Johnstone_Stone_spaces} for the special case of this construction when $A$ is a ring; it is attributed to Andr\'e Joyal in \cite[Notes on chapter V]{Johnstone_Stone_spaces}.)
Indeed, the universality of the presentation says exactly that for any distributive lattice $D$, to give a subadditive morphism $A\to D$ is equivalent to giving a lattice homomorphism $LA\to D$.

It follows that for any $A\in\OSRsa$, we have $\Rad(A)\cong \Idl(LA)$. Frames of the form $\Idl(D)$ for some $D\in \DLat$ are called \defemph{coherent} (see e.g.~\cite[Section~II.3]{Johnstone_Stone_spaces}).
The following is now immediate.
\begin{proposition}
    Let $F$ be a frame. There exists an ordered semiring $A$ with $F\cong \Rad(A)$ if and only if $F$ is coherent.
\end{proposition}

\begin{example}
    Let $B_1=(B,\bot,\top,\vee,\wedge)$ be a Boolean algebra and $B_2=(B,\bot,\top,\neg,\Delta,\wedge)$ be the corresponding Boolean ring.
    One can regard $B_2$ as a discretely ordered semiring $B_3=(B,{=},\bot,\top,\Delta,\wedge)$ (Example~\ref{ex:semiring}), and then obtain the distributive lattice $LB_3$.
    On the other hand, $B_1$ itself is also a distributive lattice.
    Since the identity function on $B$ is a subadditive morphism $B_3\to B_1$ and is in fact the universal one from $B_3$ to a distributive lattice, we have a canonical isomorphism $LB_3\cong B_1$.
\end{example}

\section{Prime ideals and spectra}
Let $A$ be an ordered semiring.
A \defemph{prime ideal} of $A$ is an ideal $I\subseteq A$ satisfying 
\begin{itemize}
    \item $1\notin I$; and 
    \item for all $x,y\in A$, if $xy\in I$ then $x\in I$ or $y\in I$.
\end{itemize}
The prime ideals of $A$ correspond to the subadditive morphisms $A\to\mathbf{2}$ (cf.~Examples~\ref{ex:morphism-dist-lat-2} and \ref{ex:morphism-ring-2}).
We denote the set of all prime ideals of $A$ by $\Spec(A)$, and call it the \defemph{spectrum} of $A$. 
We shall identify $\Spec(A)$ with the set $\OSRsa(A,\mathbf{2})$ of all subadditive morphisms $A\to\mathbf{2}$. 
We put a topology on $\Spec(A)$ as follows.
We regard $\mathbf{2}$ as the Sierpi\'nski space in which $\{1\}$ is open, endow $\mathbf{2}^A$ with the product topology, and then endow $\Spec(A)=\OSRsa(A,\mathbf{2})\subseteq \mathbf{2}^A$ with the induced topology. Explicitly, for each $x\in A$ we have an open set 
\[
D_x=\{\,(f\colon A\to \mathbf{2})\in\Spec(A)\mid f(x)=1\,\},
\]
and the family $\{\,D_x\mid x\in A\,\}$ forms a basis of the topology on $\Spec(A)$. 

Next we shall show that the topological space $\Spec(A)$ corresponds to the frame $\Rad(A)$ via the (idempotent) adjunction
\begin{equation}
\label{eqn:O-pt}
    \begin{tikzpicture}[baseline=-\the\dimexpr\fontdimen22\textfont2\relax ]
      \node(0) at (0,0) {$\mathbf{Top}$};
      \node(1) at (4,0) {$\Frm^\op$.};
      \draw [->,transform canvas={yshift=5}] (0) to node[auto, labelsize] {$\mathcal{O}$} (1); 
      \draw [<-,transform canvas={yshift=-5}] (0) to node[auto, swap,labelsize] {$\pt$} (1); 
      \path (0) to node[rotate=90] {$\vdash$} (1);
    \end{tikzpicture}
\end{equation}
Here, $\mathbf{Top}$ is the category of topological spaces and continuous maps, $\mathcal{O}$ maps each topological space $X$ to the frame $\mathcal{O}(X)$ of open sets of $X$, and $\pt$ maps each frame $F$ to the topological space $\pt(F)$, which is the set $\Frm(F,\mathbf{2})$ equipped with the topology induced from $\mathbf{2}^F$ ($\mathbf{2}$ is again regarded as the Sierpi\'nski space in which $\{1\}$ is open).
The adjunction \eqref{eqn:O-pt} restricts to an equivalence between the category of sober topological spaces and (the opposite of) the category of spatial frames;
see e.g.~\cite[Chapter~II]{Johnstone_Stone_spaces} for details.
Every coherent frame is known to be spatial \cite[Theorem~II.3.4]{Johnstone_Stone_spaces}.

\begin{proposition}
\label{prop:pt-Rad-Spec}
    Let $A$ be an ordered semiring. Then there exists a canonical homeomorphism $\pt(\Rad(A))\cong \Spec(A)$.
\end{proposition}
\begin{proof}
    The points of $\pt(\Rad(A))$ are the frame homomorphisms $\Rad(A)\to \mathbf{2}$, which corresponds to the subadditive morphisms $A\to\mathbf{2}$, i.e., to the points of $\Spec(A)$. 
    
    For each $I\in\Rad(A)$, let $U_I$ be the subset of $\pt(\Rad(A))$ defined by
    \[
    U_I=\{\,p\colon \Rad(A)\to \mathbf{2}\mid p(I)=1\,\}.
    \]
    The set of all open subsets of $\pt(\Rad(A))$ is $\{\,U_I\mid I\in\Rad(A)\,\}$.
    Under the correspondence $\pt(\Rad(A))\cong \Spec(A)$, for each $x\in A$, $D_x\subseteq \Spec(A)$ corresponds to $U_{\sqrt{\langle x\rangle}}\subseteq \pt(\Rad(A))$. Because all radical ideals are joins (in $\Rad(A)$) of radical ideals of the form $\sqrt{\langle x \rangle}$, we see that the topologies also agree.
\end{proof}

\begin{corollary}
    Let $A$ be an ordered semiring. Then there exists a canonical isomorphism of frames $\Rad(A)\cong\cat{O}(\Spec(A))$.
\end{corollary}
\begin{proof}
    Since $\Rad(A)$ is spatial, the claim follows from Proposition~\ref{prop:pt-Rad-Spec}.
\end{proof}

Modulo the isomorphism $\Rad(A)\cong\mathcal{O}(\Spec(A))$, the subadditive morphism $\sqrt{\langle - \rangle}\colon A\to \Rad(A)$ coincides with $D_{(-)}\colon A\to \mathcal{O}(\Spec(A))$.

\medbreak

We end this note with a few elementary observations on prime ideals.

First note that if $1\leq 0$ holds in an ordered semiring $A$, then $A$ has no prime ideals. 
(There exist non-trivial ordered semirings in which $1\leq 0$ holds. For example, consider the order dual of any ordered semiring obtained from Example~\ref{ex:idempotent-semiring}.)
We now show the converse.
An ideal $I$ of an ordered semiring $A$ is called a \defemph{maximal ideal} if $I\neq A$ and the only ideals containing $I$ are either $I$ or $A$.
\begin{lemma}
    Every maximal ideal in an ordered semiring is prime.
\end{lemma}
\begin{proof}
    Let $I$ be a maximal ideal of an ordered semiring $A$. Let $x,y$ be elements of $A$ such that $x\notin I$ and $y\notin I$. We shall show that $xy\notin I$.
    Since we have $\langle I\cup\{x\}\rangle=A$, there exist $z\in I$ and $u\in A$ with $1\leq z+xu$. Similarly, there exist $w\in I$ and $v\in A$ with $1\leq w+yv$. Therefore we have 
    \[
    1=1\cdot 1\leq (z+xu)(w+yv)=(zw+zyv+xuw)+xyuv\in \langle I\cup\{xy\}\rangle,
    \]
    showing that $xy\notin I$.
\end{proof}
\begin{proposition}
    Let $A$ be an ordered semiring. Then the following are equivalent.
    \begin{enumerate}[label=\emph{(\arabic*)}]
        \item $x\leq 0$ holds for any $x\in A$.
        \item $1\leq 0$ holds in $A$.
        \item $A$ has no prime ideals.
        \item $A$ has no maximal ideals.
    \end{enumerate}
\end{proposition}
\begin{proof}
    (1)$\implies$(2)$\implies$(3)$\implies$(4) is trivial, while (4)$\implies$(1) can be proved by a standard application of Zorn's lemma; note that (1) is equivalent to the statement that the smallest ideal $\langle 0\rangle =\{\,x\in A\mid x\leq 0\,\}$ of $A$ is equal to $A$.
\end{proof}

\begin{proposition}
    Let $A$ be an ordered semiring. An ideal $I\subseteq A$ is a prime ideal if and only if $I$ is a prime element of the integral quantale $\Idl(A)$, namely 
    $I\neq A$ and for any $J,K\in \Idl(A)$, $J\cdot K\subseteq I$ implies either $J\subseteq I$ or $K\subseteq I$.
\end{proposition}
\begin{proof}
    Any map $g\colon \Idl(A)\to \mathbf{2}$ yields a map $f=g\circ\langle-\rangle\colon A\to \mathbf{2}$. By Theorem~\ref{thm:Idl}, $g$ is a quantale homomorphism if and only if $f$ is a subadditive morphism. The former condition implies that $g$ preserves arbitrary joins, which is the case precisely when there exists $I\in \Idl(A)$ such that $g(J)=0$ if and only if $J\subseteq I$ for all $J\in\Idl(A)$; in this situation $f$ is the characteristic map of $I\subseteq A$.
    Now observe that $g$ is a quantale homomorphism if and only if this $I$ is a prime element of $\Idl(A)$, and $f$ is a subadditive morphism if and only if $I\subseteq A$ is a prime ideal. 
\end{proof}

\subsection*{Acknowledgements}
I thank Hisashi Aratake for providing useful comments on an earlier version of this note.

\bibliographystyle{alpha} %
\bibliography{myref} %
\end{document}